\documentclass[preprint]{elsarticle}

\usepackage{amssymb, amsmath, amsthm}
\usepackage{graphicx, epsfig, epstopdf}

\journal{xyz journal}

\begin{document}

\begin{frontmatter}

\title{The Birkhoff Diamond as Double Agent}
\author{James A. Tilley \footnotemark \textsuperscript{,} \footnotemark \textsuperscript{,} \footnotemark}

\begin{abstract}
\noindent Despite the existence of a proof of the 4-color theorem, it would seem that there is still more to learn about why any planar graph is 4-colorable. To that end, we take another look at the Birkhoff diamond and discover something new and intriguing: after an extensive search for (rare) Kempe-locked triangulations, we find a Birkhoff diamond subgraph in each one. We offer a heuristic argument as to why that result is not only reasonable but also to be expected and posit that the presence of a Birkhoff diamond is necessary to Kempe-locking. If that conjecture is true, it means that the Birkhoff diamond plays a double role in the matter of 4-colorability, simultaneously working for opposite sides of whether a given planar graph could possibly be a minimum counterexample.
\end{abstract}


\end{frontmatter}

\footnotetext[1]{address: 61 Meeting House Road, Bedford Corners, NY 10549-4238}
\footnotetext[2]{e-mail: jimtilley@optonline.net}
\footnotetext[3]{telephone: 914-242-9081}

\section*{Introduction}

\noindent Mention the Birkhoff diamond to a mathematician, particularly a graph theorist, and it brings to mind the matter of the 4-color problem. Tell her that you have some new insight about the Birkhoff diamond and she will shudder to think that you might be bold enough to step forward with a claim of a ``human'' proof (that is, without the significant aid of a computer) of the 4-color theorem. Let us be clear{---}we make no such claim. However, we do believe there is still more to understand about why any planar graph is 4-colorable and, to that end, we examine the role of the Birkhoff diamond. Needless to say, we discover something new and intriguing along the way: we find that it is plausible, perhaps highly so, that the Birkhoff diamond plays a double role in the matter of 4-colorability, simultaneously working for \textit{opposite sides} of whether a given planar graph could possibly be a minimum counterexample to the 4-color conjecture{---}hence the moniker ``double agent.'' (When we mention a minimum counterexample, we refer to the 4-color \textit{conjecture} because it makes less sense to talk about a minimum counterexample to the proved 4-color \textit{theorem}.)

Though it may appear otherwise to some, we contend that this article is not yet another in a long list of (futile) attempts to find an alternative proof of the 4-colorability of planar graphs{---}instead, it is about a way of looking at the problem to understand better what it could be about planar graphs that renders them 4-colorable. Our objective is to consider the question: ``\textit{Why} must a planar graph be 4-colorable?'' The tentative answer we find is the provocative statement: ``Because the Birkhoff diamond appears to serve two opposing masters at the same time.''

We presume a basic understanding of graph theory, but define and illustrate several less-common terms and all new ones. Any planar graph that is not a \textit{triangulation} (a graph all of whose faces are delineated by three edges) can be turned into a triangulation by inserting edges. If that triangulation can be 4-colored, so can the original graph. Hence, we focus on triangulations and their close relatives, \textit{near-triangulations}, all of whose faces \textit{except for one} are delineated by three edges. Further, we consider only graphs that are \textit{connected}{---} that is, graphs in which there is a path joining every pair of vertices. A graph is said to be \textit{$k$-connected} if it has more than $k$ vertices and remains connected \textit{whenever} fewer than $k$ vertices are deleted (see \cite{ChartrandZhang}). Because a $k$-connected graph cannot be planar if $k > 5$ (see \cite{ChartrandZhang}), we can limit our study to planar triangulations that are either 4-connected or 5-connected. An example of each is shown in figure 1. The graph on the left is the well-known icosahedron, the smallest 5-connected triangulation. We will see that the graph on the right is the smallest Kempe-locked triangulation, a term we will define in due course.

\begin{figure}[htb]
\label{figure 1}
\centering
\includegraphics*[scale=0.80, trim= 110 550 80 85] {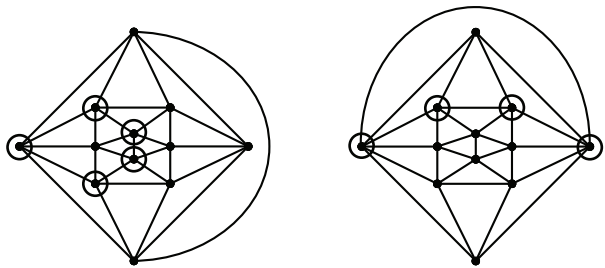}
\caption{Two planar triangulations: the one on the left is 5-connected and the one on the right is 4-connected. When the circled vertices are deleted, the graphs disconnect.}
\end{figure}

It is not necessary to consider planar triangulations that are only 3-connected because any such graph on more than four vertices must have at least one \textit{separating triangle} (a triangle with vertices of the graph both inside and outside the triangle) and a graph with a separating triangle cannot be a minimum counterexample to the 4-color conjecture. Assume such a graph $T$ is a minimum counterexample. Then, the vertex sets of both (1) the proper subgraph of $T$ consisting of the separating triangle and everything outside it and (2) the proper subgraph of $T$ consisting of the separating triangle and everything inside it can be 4-colored in such a way that the separating triangle is colored identically in each (through a permutation of colors, as necessary). Thus $T$ can be 4-colored, a contradiction.

Accepted proofs of the 4-color theorem \cite{AppelHaken,AppelHakenKoch,RobertsonSandersSeymourThomas,Steinberger} depend on two key ideas{---}unavoidability and reducibility. In this context, \textit{unavoidability} refers to a finite set of planar \textit{configurations} (near-triangulations drawn with the infinite face non-triangular) at least one of which must appear in any planar triangulation that is sufficiently connected to possibly be a counterexample to the 4-color conjecture (one that belongs to a subset of all 5-connected triangulations in which the removal of any 5-cycle disconnects the graph into two components, one of which is a single vertex{---}for example, the left panel of figure 1). In other words, the set of configurations is unavoidable. A \textit{reducible} configuration is one whose presence in a planar triangulation renders that triangulation 4-colorable if the near-triangulation with the reducible configuration removed is 4-colorable. Thus, a reducible configuration cannot appear in a minimum counterexample to the 4-color conjecture. Consequently, if one can create an unavoidable set of reducible configurations, there can be no minimum counterexample and the 4-color conjecture is proved. The Birkhoff diamond, depicted in figure 2, is named for the mathematician who proved that it is a reducible configuration \cite{Birkhoff}. It is the smallest configuration appearing in every one of the unavoidable sets of reducible configurations used to prove the 4-color theorem \cite{AppelHaken,AppelHakenKoch,RobertsonSandersSeymourThomas,Steinberger}.

\begin{figure}[htb]
\label{figure 2}
\centering
\includegraphics*[scale=0.80, trim= 90 600 75 100] {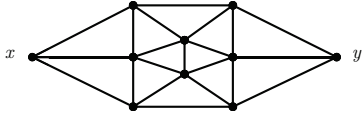}
\caption{The Birkhoff diamond with endpoints $x$ and $y$.}
\end{figure}

We have just described the known first role of the Birkhoff diamond in the matter of 4-colorability of planar graphs: its \textit{mere presence} in a planar triangulation serves to disqualify that graph from being a minimum counterexample to the 4-color conjecture. We now turn to the posited second role, the primary subject of this article: its presence in a planar triangulation is required for that graph to possibly be a minimum counterexample. If this second role could be verified, then the 4-color theorem would follow trivially. We therefore suspect that any proof of the second role would be nontrivial and require superhuman effort. Nevertheless, because our extensive testing has produced strongly suggestive results, we claim that there is new and valuable insight to be gained from the conjectured second role, even absent a proof of its validity.

To arrive at the supposition regarding the second role of the Birkhoff diamond, we first need to introduce the idea of Kempe-locking and for that we need to know what a Kempe chain is. Those are the subjects of the next two sections. With those concepts in hand, we are able to demonstrate that a minimum counterexample to the 4-color conjecture must be Kempe-locked with respect each of its edges. This is a highly restrictive condition, one that gets harder and harder to meet as the size of a triangulation increases. (Because the 4-color conjecture has been proved, we know that the condition is actually impossible to satisfy.) The results of an extensive search for Kempe-locked triangulations lead to the plausible conjecture that a planar triangulation cannot be Kempe-locked with respect to an edge unless the vertices serving as endpoints of that edge also serve as the ``endpoints'' of a Birkhoff diamond ($x$ and $y$ in figure 2).

\section*{Kempe chains}

\noindent An important tool in graph coloring is the \textit{Kempe chain}, named after the British mathematician whose famous attempt at proving the 4-color conjecture failed \cite{Kempe}. We deal only with \textit{proper} vertex-colorings of a graph, those in which \textit{adjacent} vertices (those joined by an edge) must have different colors. Following the standard convention, we use the integers 1, 2, 3, 4, \dots  to indicate distinct colors. In a given coloring of a graph $G$, a Kempe chain is a maximal, connected, induced subgraph of $G$ whose vertices use only two colors, let us say colors $i$ and $j$. (An \textit{induced} subgraph $F$ of a graph $G$ is one in which all edges in $G$ that join vertices in the vertex set of $F$ are also edges in $F$.) An $i$-$j$ Kempe chain is ``maximal'' in the sense that every vertex adjacent to, but not in, the chain uses a color other than $i$ or $j$. A \textit{short $i$-$j$ Kempe chain} consists of a single vertex and uses color $i$ or color $j$. All vertices adjacent to such a single-vertex Kempe chain use colors other than $i$ or $j$. Kempe chains are illustrated in figure 3.

\begin{figure}[htb]
\label{figure 3}
\centering
\includegraphics*[scale=0.80, trim= 90 510 80 90] {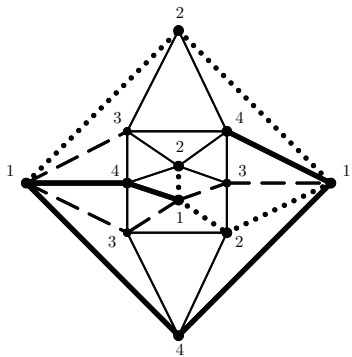}
\caption{In this 4-coloring of a near-triangulation, three Kempe chains have been highlighted: the 1-2 chain by dotted lines, the 1-3 chain by dashed lines, and the 1-4 chain by solid lines.}
\end{figure}

Kempe chains are particularly useful for ``navigating'' among a subset (possibly the whole set) of the distinct 4-colorings of $G$ because interchanging colors on a chain{---}that is, interchanging the color labels $i$ and $j$ for all vertices constituting an $i$-$j$ chain{---}leaves $G$ properly colored. Interchanging colors on an $i$-$j$ Kempe chain does not result in a distinctly different coloring of $G$ if there is only one $i$-$j$ chain in $G$. Any proper vertex-coloring of a graph partitions the vertex set of the graph into \textit{color classes} and when there is only one $i$-$j$ Kempe chain, a color interchange on that chain does not alter the partitioning of the vertex set into color classes.

\section*{Kempe-locked triangulations}

\noindent Kempe-locking is a property of a planar triangulation $T$ with respect to one of its edges $xy$. Let $G_{xy}$ denote the near-triangulation that results when the edge $xy$ is deleted from $T$: $T$ is said to be \textit{Kempe-locked with respect to the edge $xy$} if, in \textit{every} 4-coloring of $G_{xy}$ in which the colors of $x$ and $y$ are the same, there are precisely \textit{three} Kempe chains that include both $x$ and $y$. Thus, if $T$ is Kempe-locked with respect to the edge $xy$, then given any 4-coloring of $G_{xy}$ in which $x$ and $y$ are both colored the same, without loss of generality color 1, there must be 1-2, 1-3, and 1-4 Kempe chains including both $x$ and $y$. Interchanging colors on any of those chains leaves $x$ and $y$ sharing the same color. Moreover, because $T$ is Kempe-locked with respect to $xy$, interchanging colors on any Kempe chain involving neither $x$ nor  $y$ results in a 4-coloring in which there are still 1-2, 1-3, and 1-4 Kempe chains that include both $x$ and $y$. If interchanging colors on Kempe chains is the only method of recoloring available, $G_{xy}$, once in the state in which $x$ and $y$ have the same color, is ``locked into'' that state.

\begin{figure}[htb]
\label{figure 4}
\centering
\includegraphics*[scale=0.80, trim= 90 505 80 90] {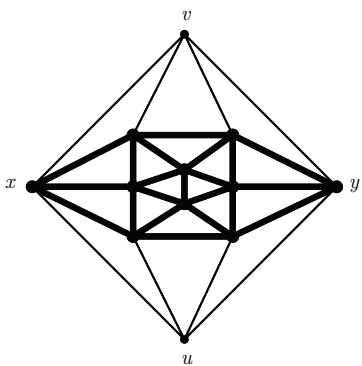}
\caption{The near-triangulation $G_{xy}$ of order 12 obtained when the edge that had joined $x$ and $y$ is deleted from the triangulation shown in the right panel of figure 1. The Birkhoff diamond from figure 2 is shown highlighted as a subgraph of $G_{xy}$. A 4-coloring of $G_{xy}$ in which $x$ and $y$ are colored the same is given in figure 3.}
\end{figure}

Figure 4 shows the near-triangulation $G_{xy}$ that results from deleting the edge $xy$  in the triangulation $T$ on 12 vertices (we say that $T$ is of \textit{order} 12) illustrated in the right panel of figure 1. We adopt a convention, illustrated in figure 4, of drawing $G_{xy}$ with the 4-face as the exterior (infinite) face and with $x$ and $y$ denoting the leftmost and rightmost vertices, respectively, on the boundary of that face. We label the boundary 4-cycle $uxvy$ with $u$ the bottommost vertex and $v$ the topmost vertex. The $G_{xy}$ depicted in figure 4 has as a subgraph the Birkhoff diamond from figure 2{---}it is highlighted. Figure 3 gives a proper 4-coloring of this $G_{xy}$ in which $x$ and $y$ are both colored 1. From figure 3, we note that there are 1-2 and 1-4 Kempe chains (each of which includes part of the the boundary of $G_{xy}$ as drawn) that include both $x$ and $y$ and a 1-3 Kempe chain including both $x$ and $y$ that snakes all the way through the interior of $G_{xy}$. It is easily verified that there are only five other distinct 4-colorings of this $G_{xy}$ with both $x$ and $y$ colored 1. They are readily found once $v$ is colored 2 and the upper boundary of the Birkhoff diamond is colored 1-3-4-1 as in figure 3, no loss of generality with any of these color choices. In each of those five distinct additional 4-colorings of $G_{xy}$, there are 1-2, 1-3, and 1-4 Kempe chains including both $x$ and $y$. Thus, the planar triangulation $T$ from which this $G_{xy}$ is derived (the right panel of figure 1) is Kempe-locked with respect to the edge $xy$. It features a Birkhoff diamond subgraph with endpoints $x$ and $y$.

We now show that a minimum counterexample to the $4$-color conjecture must be Kempe-locked with respect to each of its edges. Let a planar triangulation $T$ be a minimum counterexample to the 4-color conjecture. Consider an arbitrary edge $xy$ and coalesce $x$ and $y$ into $w$ so that all the edges formerly incident to $x$ and $y$ now become incident to $w$. This so-called \textit{edge contraction} yields a new triangulation $T^{\prime}$ with one fewer vertex than $T$. Because $T$ is assumed to be a minimum counterexample, $T^{\prime}$ can be 4-colored. Then, when $w$ is split apart into the original $x$ and $y$, but without replacing the edge $xy$, we obtain a near-triangulation $G_{xy}$ in which the colors of $x$ and $y$ are the same. Now suppose that $T$ is not Kempe-locked with respect to the edge $xy$. Then there must be a 4-coloring of $G_{xy}$ in which $x$ and $y$ are colored the same \textit{and} in which there is a Kempe chain that includes $x$ but does not include $y$. Interchanging colors on such a chain results in a 4-colored $G_{xy}$ with the color of $x$ not the same as the color of $y$. In this 4-coloring, the edge $xy$ can be inserted to yield a 4-coloring for $T$, in contradiction to the assumption that $T$ is a minimum counterexample.

\section*{Fundamental Kempe-locking configurations}

\noindent At last we are in a position to identify configurations that are critical to Kempe-locking. Using the same notation as previously, let $T$ be a planar triangulation and let $G_{xy}$ be the near-triangulation derived from it by deleting the edge $xy$. Finally, with $uxvy$ the 4-cycle delineating the infinite 4-face of $G_{xy}$, let $K_{xy}$ be the near-triangulation that results when the bottommost and topmost vertices $u$ and $v$, respectively, and the edges incident to them, are deleted. We refer to $K_{xy}$ as the \textit{Kempe-locking configuration for $T$ with respect to the edge $xy$}. If $T$ is of order $n$, then $K_{xy}$ is of order $n - 2$. Are there certain Kempe-locking configurations that are more fundamental than others? Indeed. If a Kempe-clocking configuration $K_{xy}$ has no subgraph $K^{\prime}_{xy}$ that is the Kempe-locking configuration with respect to the edge $xy$ for some planar triangulation $T^{\prime}$ with an edge $xy$, then we say that $K_{xy}$ is a \textit{fundamental} Kempe-locking configuration.

From the previous section we see that the Birkhoff diamond is a Kempe-locking configuration and, as we will learn shortly, it turns out to be a fundamental Kempe-locking configuration. Refer to figure 5 for an example of a Kempe-locking configuration that is not fundamental. Are there any Kempe-locking configurations other than the Birkhoff diamond that are fundamental? Investigating this question is an important step toward determining whether a minimum counterexample can exist. We have shown that a minimum counterexample must be Kempe-locked with respect to each of its edges. Hence, the vertices serving as the endpoints of any given edge in a minimum counterexample must also serve as the endpoints of a fundamental Kempe-locking configuration, which may or may not be a proper subgraph of a Kempe-locking configuration. Because a triangulation of order $n$ has $3n - 6$ edges (see \cite{ChartrandZhang}), it would seem that finding a minimum counterexample of order $n$ becomes less and less likely as $n$ increases. Likewise, so it would seem that finding a fundamental Kempe-locking configuration becomes less and less likely the higher its order.

\begin{figure}[htb]
\label{figure 5}
\centering
\includegraphics*[scale=0.80, trim= 95 500 70 95] {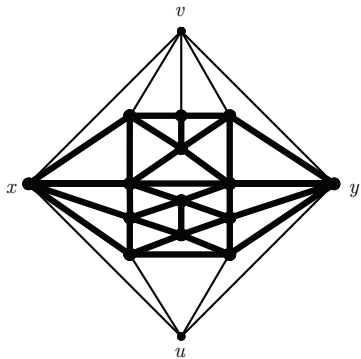}
\caption{A near-triangulation $G_{xy}$ that becomes a Kempe-locked triangulation when vertices $x$ and $y$ on the infinite 4-face are joined. The highlighted Kempe-locking configuration $K_{xy}$ is \textit{not} fundamental because it has a Birkhoff diamond subgraph with endpoints $x$ and $y$.}
\end{figure}

Consider a coloring of $G_{xy}$ of order $n$ in which $x$ and $y$ are both colored $k$. If there is a 2-color path between $v$ and $u$ that uses colors $i,j \neq k$, then there cannot be three Kempe chains that include both $x$ and $y$. A fundamental Kempe-locking configuration $K_{xy}$ must be able to ``block'' the passage of such 2-color paths from $v$ to $u$ in all 4-colorings of $G_{xy}$ in which $x$ and $y$ have the same color. There are two ways in which this so-called blocking can occur: the fundamental Kempe-locking configuration can (1) prevent the transmission of the 2-color path from the configuration's top boundary to its bottom boundary on the way from $v$ to $u$ or (2) ultimately force $u$ to take a color other than $i$ or $j$. Refer to figure 3 to see how the Birkhoff diamond configuration is able to block the 2-4 path by way of (1) and how it is able to block the 2-3 path by way of (2). Either (1) or (2) will occur whenever there is Kempe-locking.

As the order of $G_{xy}$ increases it becomes increasingly unlikely to find a fundamental Kempe-locking configuration $K_{xy}$ because the number of 4-colorings of $G_{xy}$ with $x$ and $y$ colored the same grows rapidly and the probability that there will be no 4-coloring at all with a relevant 2-color path between $v$ and $u$ correspondingly diminishes rapidly. However, since we shall learn that the Birkhoff diamond is a fundamental Kempe-locking configuration, we should expect to see it appear as a subgraph of larger and larger Kempe-locked triangulations due to its ability to block 2-color paths. Indeed, that is what we shall discover.

The approach that we adopted in the search for \textit{fundamental} Kempe-locking configurations is analogous to that of experimental physicists in their search for a new \textit{elementary} particle with specified properties. Physicists confine their explorations to collision events involving total energy in an interval sufficient to bring the sought-after particle into being. Similarly, we explore graphs of orders in which \textit{fundamental} Kempe-locking configurations would be expected to be found. The best chance to discover fundamental Kempe-locking configurations would seem to occur when there is a small number of distinct 4-colorings of a near-triangulation $G_{xy}$ with $x$ and $y$ the same color, thus admitting the possibility that \textit{every one} of those few colorings will feature three Kempe chains that include both $x$ and $y$. We are naturally led to consider low-order triangulations that are either 4-connected or 5-connected.

To assure that we did not miss any low-order triangulations, we generated the full set of 8,044 isomorphism classes for 4-connected triangulations of orders 6-15 and the full set of 9,733 isomorphism classes for 5-connected triangulations of orders 12-24. (See \cite{BrinkmannMcKay1} and \cite{BrinkmannMcKay2}.) We then tested \textit{every edge} in a triangulation from each of the isomorphism classes to determine if there are any 4-colorings with $x$ and $y$ colored the same but with fewer than three Kempe chains that include both $x$ and $y$. It turns out to be rare for this not to be the case. The only Kempe-locked triangulations we encountered were 4-connected; there were none at all among 5-connected triangulations. There are no Kempe-locked triangulations of order less than 12, a single 4-connected Kempe-locked triangulation of order 12 (the right panel of figure 1), none of order 13, a single one of order 14, and a single one of order 15. Those three Kempe-locked triangulations all feature a Birkhoff diamond configuration with $x$ and $y$ as endpoints. Each is Kempe-locked with respect to only a single edge.

In an expanded search for fundamental Kempe-locking configurations, we examined 4-connected triangulations of orders 16-20. Because the number of isomorphism classes grows rapidly with increasing order (from 30,926 at order 16 to 24,649,284 at order 20{---}refer to \cite{BrinkmannMcKay2}) and because the number of edges in a triangulation increases with increasing order, we soon ran into computation-time limitations imposed by our laptop computer. After deciding to limit aggregate computer execution time to weeks instead of months, we proceeded in the expanded search by examining all 30,926 isomorphism classes of order 16 and all 158,428  isomorphism classes of order 17, but only 100,000 randomly generated non-isomorphic triangulations for each order from 18 through 20.

For orders 16 and 17, we discovered \textit{eight}  and \textit{fourteen} non-isomorphic triangulations, respectively, that are Kempe-locked, all with respect to a single edge, call it $xy$ in each case. Each of those 4-connected Kempe-locked triangulations features a Birkhoff diamond with $x$ and $y$ as endpoints.  Figure 5 shows the $G_{xy}$ derived from one of the Kempe-locked triangulations of order 16. In the random samples of 100,000 triangulations each for orders 18-20, we discovered additional non-isomorphic Kempe-locked triangulations: \textit{ten} of order 18, \textit{eight} of order 19, and \textit{five} of order 20, all locked with respect to a single edge $xy$ and all featuring a Birkhoff diamond configuration with $x$ and $y$ as endpoints. It is an open question whether there are any triangulations that are at least 4-connected and Kempe-locked with respect to more than a single edge.

Let us take stock of the results of our search for fundamental Kempe-locking configurations. We found no 5-connected Kempe-locked triangulations of orders 12-24. The only fundamental Kempe-locking configuration we found is the Birkhoff diamond of order 10. We have shown that there are no fundamental Kempe-locking configurations of orders 9 or less or between 11 and 15, inclusive, and a sample of 100,000 non-isomorphic 4-connected triangulations each for orders 18-20 turned up no fundamental Kempe-locking configurations of orders 16-18. We conclude that the ``experimental'' case is strong that the Birkhoff diamond is the only fundamental Kempe-locking configuration.

\section*{Conclusion}

\noindent The results of our search for fundamental Kempe-locking configurations point to a reasonable conjecture that the Birkhoff diamond is the only one. An equivalent way to state this conjecture is that the existence of a Birkhoff diamond configuration with endpoints $x$ and $y$ in a planar triangulation $T$ with an edge $xy$ is a necessary, but not sufficient, condition for $T$ to be Kempe-locked with respect to the edge $xy$.  This is the conjectured second role of the Birkhoff diamond{---}as the ``agent'' for the existence of a minimum counterexample to the 4-color conjecture.

A minimum counterexample to the 4-color conjecture must be Kempe-locked with respect to each of its edges and if the conjecture regarding the Birkhoff diamond is true, then a minimum counterexample to the 4-color conjecture must have at least as many Birkhoff diamond subgraphs present as there are edges, any one of which would render the minimum counterexample 4-colorable and hence not a minimum counterexample. Moreover, the existence of that many Kempe-locking Birkhoff diamonds in a planar triangulation is easily seen to be impossible because the degrees of $x$ and $y$ in an edge $xy$ that is Kempe-locked by a Birkhoff diamond must be at least 6: having as many such Birkhoff diamonds as there are edges would mean that the triangulation is \textit{minimum degree} 6 and thus nonplanar (see \cite{ChartrandZhang}).  (In the icosahedron there are as many Birkhoff diamonds as edges, but none is Kempe-locking.) Indeed, such an every-edge-Kempe-locked triangulation cannot be constructed because the signature feature of a Birkhoff diamond is its central diamond consisting of four vertices of degree 5{---}thus, none of the five edges of that central diamond can support a Birkhoff-diamond Kempe-locking configuration.

For a planar triangulation to be a minimum counterexample, it appears necessary for it to have a Birkhoff diamond subgraph, but if it does, then it cannot be a minimum counterexample. As remarkable as the Birkhoff diamond is, only its imaginary ``quantum'' version can be both present and absent at the same time. As stated in the introduction, it will likely be very hard to prove that the Birkhoff diamond is the only fundamental Kempe-locking configuration, a situation that is unfortunately all too common in mathematics{---}a conjecture that is easy to formulate, has plenty of supporting evidence (namely, despite efforts to find one, there is no known counterexample), is highly difficult to prove, and is possibly false. Nevertheless, the conjecture posited in this article, if true, has the satisfying feature that the Birkhoff diamond configuration alone would explain why any planar graph is 4-colorable. That makes the conjecture worthy of further study.


\end{document}